# The distribution of model averaging estimators and an impossibility result regarding its estimation

**Benedikt M. Pötscher**[1]

*University of Vienna*

**Abstract:** The finite-sample as well as the asymptotic distribution of Leung and Barron's (2006) model averaging estimator are derived in the context of a linear regression model. An impossibility result regarding the estimation of the finite-sample distribution of the model averaging estimator is obtained.

## 1. Introduction

Model averaging or model mixing estimators have received increased interest in recent years; see, e.g., Yang [18, 19, 20], Magnus [13], Leung and Barron [12], and the references therein. [For a discussion of model averaging from a Bayesian perspective see Hoeting et al. [4].] The main idea behind this class of estimators is that averaging estimators obtained from different models should have the potential to achieve better overall risk performance when compared to a strategy that only uses the estimator obtained from one model. As a consequence, the above mentioned literature concentrates on studying the risk properties of model averaging estimators and on associated oracle inequalities. In this paper we derive the finite-sample as well as the asymptotic distribution (under fixed as well as under moving parameters) of the model averaging estimator studied in [12]; for the sake of simplicity we concentrate on the special case when only two candidate models are considered. Not too surprisingly, it turns out that the finite-sample distribution (after centering and scaling) depends on unknown parameters, and thus cannot be directly used for inferential purposes. As a consequence, one may be interested in estimators of this distribution, e.g., for purposes of conducting inference. We establish an impossibility result by showing that any estimator of the finite-sample distribution of the model averaging estimator is necessarily "bad" in a sense made precise in Section 4. While we concentrate on Leung and Barron's [12] estimator (in the context of only two candidate models) as a prototypical example of a model averaging estimator in this paper, similar results will typically hold for other model averaging estimators (and more than two candidate models) as well.

We note that results on distributional properties of post-model-selection estimators that parallel the development in the present paper have been obtained in [5, 6, 7, 9, 10, 14, 15, 16, 17]. See also Leeb and Pötscher [11] for impossibility results pertaining to shrinkage-type estimators like the Lasso or Stein's estimator. An easily accessible exposition of the issues discussed in the just mentioned literature can be found in Leeb and Pötscher [8].

---

[1]Department of Statistics, University of Vienna, e-mail: benedikt.poetscher@univie.ac.at
*AMS 2000 subject classifications:* primary 62F10, 62F12; secondary 62E15, 62J05, 62J07.
*Keywords and phrases:* model mixing, model aggregation, combination of estimators, model selection, finite sample distribution, asymptotic distribution, estimation of distribution.

The only other paper we are aware of that considers distributional properties of model averaging estimators is Hjort and Claeskens [3]. Hjort and Claeskens [3] provide a result (Theorem 4.1) that says that – under some regularity conditions – the asymptotic distribution of a model averaging estimation scheme is the distribution of the same estimation scheme applied to the limiting experiment (which is a multivariate normal estimation problem). This result is an immediate consequence of the continuous mapping theorem, and furthermore becomes vacuous if the estimation problem one starts with is already a Gaussian problem (as is the case in the present paper).

## 2. The model averaging estimator and its finite-sample distribution

Consider the linear regression model

$$Y = X\beta + u$$

where $Y$ is $n \times 1$ and where the $n \times k$ non-stochastic design matrix $X$ has full column rank $k$, implying $n \geq k$. Furthermore, $u$ is normally distributed $N(0, \sigma^2 I_n)$, $0 < \sigma^2 < \infty$. Although not explicitly shown in the notation, the elements of $Y$, $X$, and $u$ may depend on sample size $n$. [In fact, the random variables $Y$ and $u$ may be defined on a sample space that varies with $n$.] Let $\mathbb{P}_{n,\beta,\sigma}$ denote the probability measure on $\mathbb{R}^n$ induced by $Y$, and let $\mathbb{E}_{n,\beta,\sigma}$ denote the corresponding expectation operator. As in [12], we also assume that $\sigma^2$ is known (and thus is fixed). [Results for the case of unknown $\sigma^2$ that parallel the results in the present paper can be obtained if $\sigma^2$ is replaced by the residual variance estimator derived from the unrestricted model. The key to such results is the observation that this variance estimator is independent of the least squares estimator for $\beta$. The same idea has been used in [7] to derive distributional properties of post-model-selection estimators in the unknown variance case from the known variance case. For brevity we do not give any details on the unknown variance case in this paper.] Suppose further that $k > 1$, and that $X$ and $\beta$ are commensurably partitioned as

$$X = [X_1 : X_2]$$

and $\beta = [\beta_1', \beta_2']'$ where $X_i$ has dimension $k_i \geq 1$. Let the restricted model be defined as $M_R = \{\beta \in \mathbb{R}^k : \beta_2 = 0\}$ and let $M_U = \mathbb{R}^k$ denote the unrestricted model. Let $\hat{\beta}(R)$ denote the restricted least squares estimator, i.e., the $k \times 1$ vector given by

$$\hat{\beta}(R) = \begin{bmatrix} (X_1'X_1)^{-1}X_1'Y \\ 0_{k_2 \times 1} \end{bmatrix},$$

and let $\hat{\beta}(U) = (X'X)^{-1}X'Y$ denote the unrestricted least squares estimator. Leung and Barron [12] consider model averaging estimators in a linear regression framework allowing for more than two candidate models. Specializing their estimator to the present situation gives

$$\tilde{\beta} = \hat{\lambda}\hat{\beta}(R) + (1 - \hat{\lambda})\hat{\beta}(U) \tag{1}$$

where the weights are given by

$$\hat{\lambda} = [\exp(-\alpha\hat{r}(R)/\sigma^2) + \exp(-\alpha\hat{r}(U)/\sigma^2)]^{-1} \exp(-\alpha\hat{r}(R)/\sigma^2).$$

Here $\alpha > 0$ is a tuning parameter (note that Leung and Barron's tuning parameter corresponds to $2\alpha$) and

$$\hat{r}(R) = Y'Y - \hat{\beta}(R)'X'X\hat{\beta}(R) + \sigma^2(2k_1 - n)$$

and

$$\hat{r}(U) = Y'Y - \hat{\beta}(U)'X'X\hat{\beta}(U) + \sigma^2(2k - n).$$

For later use we note that

$$\begin{aligned} \hat{\lambda} &= [1 + \exp(-2\alpha k_2)\exp(-\alpha(\hat{\beta}(R)'X'X\hat{\beta}(R) - \hat{\beta}(U)'X'X\hat{\beta}(U))/\sigma^2)]^{-1} \\ &= [1 + \exp(-2\alpha k_2)\exp(\alpha \left\| X\hat{\beta}(R) - X\hat{\beta}(U) \right\|^2 /\sigma^2)]^{-1} \end{aligned} \tag{2}$$

where $\|x\|$ denotes the Euclidean norm of a vector $x$, i.e., $\|x\| = (x'x)^{1/2}$. Leung and Barron [12] establish an oracle inequality for the risk $\mathbb{E}_{n,\beta,\sigma}(\|X(\tilde{\beta} - \beta)\|^2)$ and show that the model averaging estimator performs favourably in terms of this risk. As noted in the introduction, in the present paper we consider distributional properties of this estimator. Before we now turn to the finite-sample distribution of the model averaging estimator we introduce some notation: For a symmetric positive definite matrix $A$ the unique symmetric positive definite root is denoted by $A^{1/2}$. The largest (smallest) eigenvalue of a matrix $A$ is denoted by $\lambda_{\max}(A)$ ($\lambda_{\min}(A)$). Furthermore, $P_R$ and $P_U$ denote the projections on the column space of $X_1$ and of $X$, respectively.

**Proposition 1.** *The finite-sample distribution of $\sqrt{n}(\tilde{\beta} - \beta)$ is given by the distribution of*

$$\begin{aligned} & B_n\sqrt{n}\beta_2 + C_n\sqrt{n}Z_1 + \\ & \left[1 + \exp(2\alpha k_2)\exp\left(-\alpha \left\| Z_2 + (X_2'(I - P_R)X_2)^{1/2}\beta_2 \right\|^2 /\sigma^2\right)\right]^{-1} \times \\ & \{D_n\sqrt{n}Z_2 - B_n\sqrt{n}\beta_2\} \end{aligned} \tag{3}$$

*which can also be written as*

$$\begin{aligned} & C_n\sqrt{n}Z_1 + D_n\sqrt{n}Z_2 - \\ & \left[1 + \exp(-2\alpha k_2)\exp\left(\alpha \left\| Z_2 + (X_2'(I - P_R)X_2)^{1/2}\beta_2 \right\|^2 /\sigma^2\right)\right]^{-1} \times \\ & \{D_n\sqrt{n}Z_2 - B_n\sqrt{n}\beta_2\}. \end{aligned} \tag{4}$$

*Here*

$$B_n = \begin{bmatrix} (X_1'X_1)^{-1}X_1'X_2 \\ -I_{k_2} \end{bmatrix}, \qquad C_n = \begin{bmatrix} (X_1'X_1)^{-1/2} \\ 0_{k_2 \times k_1} \end{bmatrix},$$

$$D_n = \begin{bmatrix} -(X_1'X_1)^{-1}X_1'X_2(X_2'(I - P_R)X_2)^{-1/2} \\ (X_2'(I - P_R)X_2)^{-1/2} \end{bmatrix},$$

*and $Z_1$ and $Z_2$ are independent, $Z_1 \sim N(0, \sigma^2 I_{k_1})$, and $Z_2 \sim N(0, \sigma^2 I_{k_2})$.*

*Proof.* Observe that

$$\tilde{\beta} = \hat{\beta}(R) + (1 - \hat{\lambda})(\hat{\beta}(U) - \hat{\beta}(R)) = \hat{\beta}(R) + (1 - \hat{\lambda})(X'X)^{-1}X'(P_U - P_R)Y$$

with $P_R = X_1(X_1'X_1)^{-1}X_1'$ and $P_U = X(X'X)^{-1}X'$. Diagonalize the projection matrix $P_U - P_R$ as

$$P_U - P_R = \mathcal{U}\Delta\mathcal{U}'$$

where the orthogonal $n \times n$ matrix $\mathcal{U}$ is given by

$$\mathcal{U} = [\mathcal{U}_1, \mathcal{U}_2, \mathcal{U}_3] = \left[X_1(X_1'X_1)^{-1/2} : (I - P_R)X_2(X_2'(I - P_R)X_2)^{-1/2} : \mathcal{U}_3\right]$$

with $\mathcal{U}_3$ representing an $n \times (n-k)$ matrix whose columns form an orthonormal basis of the orthogonal complement of the space spanned by the columns of $X$. The $n \times n$ matrix $\Delta$ is diagonal with the first $k_1$ as well as the last $n-k$ diagonal elements equal to zero, and the remaining $k_2$ diagonal elements being equal to 1. Furthermore, set $V = \mathcal{U}'Y$ which is distributed $N(\mathcal{U}'X\beta, \sigma^2 I_n)$. Then

$$\left\| X\hat{\beta}(U) - X\hat{\beta}(R) \right\|^2 = \|(P_U - P_R)Y\|^2 = \|\Delta V\|^2 = \|V_2\|^2$$

where $V_2$ is taken from the partition of $V' = (V_1', V_2', V_3')'$ into subvectors of dimensions $k_1$, $k_2$, and $n-k$, respectively. Note that $V_2$ is distributed $N((X_2'(I - P_R)X_2)^{1/2}\beta_2, \sigma^2 I_{k_2})$. Hence, in view of (2) we have that $(1-\hat{\lambda})(\hat{\beta}(U) - \hat{\beta}(R))$ is equal to

$$\begin{aligned}
&\left[1 + \exp(2\alpha k_2)\exp\left(-\alpha \|V_2\|^2/\sigma^2\right)\right]^{-1}(X'X)^{-1}X'\mathcal{U}\Delta V \\
&= \left[1 + \exp(2\alpha k_2)\exp\left(-\alpha \|V_2\|^2/\sigma^2\right)\right]^{-1}(X'X)^{-1}\begin{bmatrix} 0_{k_1\times 1} \\ X_2'\mathcal{U}_2 V_2 \end{bmatrix} \\
&= \left[1 + \exp(2\alpha k_2)\exp\left(-\alpha \|V_2\|^2/\sigma^2\right)\right]^{-1} D_n V_2.
\end{aligned}$$

Furthermore,

$$\begin{aligned}
\hat{\beta}(R) &= (X'X)^{-1}X'P_R Y \\
&= (X'X)^{-1}X'P_R\mathcal{U}V \\
&= (X'X)^{-1}X'X_1(X_1'X_1)^{-1/2}V_1 \\
&= \begin{bmatrix} (X_1'X_1)^{-1/2}V_1 \\ 0_{k_2\times 1} \end{bmatrix} = C_n V_1
\end{aligned}$$

with $V_1$ distributed $N((X_1'X_1)^{-1/2}X_1'X\beta, \sigma^2 I_{k_1})$. Hence, the finite sample distribution of $\tilde{\beta}$ is the distribution of

$$C_n V_1 + \left[1 + \exp(2\alpha k_2)\exp\left(-\alpha \|V_2\|^2/\sigma^2\right)\right]^{-1} D_n V_2 \tag{5}$$

where $V_1$ and $V_2$ are independent normally distributed with parameters given above. Defining $Z_i$ as the centered versions of $V_i$, subtracting $\beta$, and scaling by $\sqrt{n}$ then delivers the result. □

**Remark 2.** (i) The first two terms in (3) represent the distribution of $\sqrt{n}(\hat{\beta}(R) - \beta)$, whereas the third term represents the distribution of $(1-\hat{\lambda})\sqrt{n}(\hat{\beta}(U) - \hat{\beta}(R))$. In (4), the first two terms represent the distribution of $\sqrt{n}(\hat{\beta}(U) - \beta)$, whereas the third term represents the distribution of $-\hat{\lambda}\sqrt{n}(\hat{\beta}(U) - \hat{\beta}(R))$.

(ii) If $\beta_2 = 0$ then (3) can be rewritten as

$$C_n\sqrt{n}Z_1 + \|Z_2\|\left[1 + \exp(2\alpha k_2)\exp\left(-\alpha \|Z_2\|^2/\sigma^2\right)\right]^{-1} D_n\sqrt{n}(Z_2/\|Z_2\|)$$

showing that this term has the same distribution as

$$C_n\sqrt{n}Z_1 + \sqrt{\chi^2}[1+\exp(2\alpha k_2)\exp\left(-\alpha\chi^2/\sigma^2\right)]^{-1}D_n\sqrt{n}U$$

where $\chi^2$ is distributed as a $\chi^2$ with $k_2$ degrees of freedom, $U = Z_2/\left\|Z_2\right\|$ is uniformly distributed on the unit sphere in $\mathbb{R}^{k_2}$, and $Z_1$, $\chi^2$, and $U$ are mutually independent.

**Theorem 3.** *The finite-sample distribution of* $\sqrt{n}(\tilde{\beta}-\beta)$ *possesses a density* $f_{n,\beta,\sigma}$ *given by*

$$\begin{aligned}
f_{n,\beta,\sigma}(t) &= (2\pi\sigma^2)^{-k/2}[\det(X'X/n)]^{1/2} \\
&\times \exp\left(-(2\sigma^2)^{-1}\left\|n^{-1/2}(X_1'X_1)^{1/2}t_1 + n^{-1/2}(X_1'X_1)^{-1/2}X_1'X_2t_2\right\|^2\right) \\
&\times \left[1+\exp\left(-\alpha\sigma^{-2}g\left(\left\|n^{-1/2}D_{n2}^{-1}(t_2+n^{1/2}\beta_2)\right\|\right)^2 + 2\alpha k_2\right)\right]^{k_2} \\
&\times \left\{1+2\alpha\sigma^{-2}g\left(\left\|n^{-1/2}D_{n2}^{-1}(t_2+n^{1/2}\beta_2)\right\|\right)^2\right. \\
&\times \left.\left[1+\exp\left(\alpha\sigma^{-2}g\left(\left\|n^{-1/2}D_{n2}^{-1}(t_2+n^{1/2}\beta_2)\right\|\right)^2 - 2\alpha k_2\right)\right]^{-1}\right\}^{-1} \\
&\times \exp\left(-(2\sigma^2)^{-1}\left\|g\left(\left\|n^{-1/2}D_{n2}^{-1}(t_2+n^{1/2}\beta_2)\right\|\right)\right.\right. \\
&\times \left.\left.\left\|n^{-1/2}D_{n2}^{-1}(t_2+n^{1/2}\beta_2)\right\|^{-1} n^{-1/2}D_{n2}^{-1}(t_2+n^{1/2}\beta_2) - D_{n2}^{-1}\beta_2)\right\|^2\right),
\end{aligned} \tag{6}$$

*where $t$ is partitioned as $(t_1', t_2')'$ with $t_1$ being a $k_1\times 1$ vector. Furthermore, $D_{n2} = (X_2'(I-P_R)X_2)^{-1/2}$, and $g$ is as defined in the Appendix (with $a = \exp(2\alpha k_2)$ and $b = \alpha^{-1}\sigma^2$).*

*Proof.* By (5) we have that the finite-sample distribution of $\sqrt{n}(\tilde{\beta}-\beta)$ is the distribution of

$$-\sqrt{n}\beta + \sqrt{n}[C_n : D_n][V_1' : V_3']'$$

where

$$V_3 = \left[1+\exp(2\alpha k_2)\exp\left(-\alpha\left\|V_2\right\|^2/\sigma^2\right)\right]^{-1}V_2.$$

By Lemmata 15 and 16 in the Appendix it follows that $V_3$ possesses the density

$$\begin{aligned}
\psi(v_3) &= (2\pi\sigma^2)^{-k_2/2}\left[1+\exp\left(-\alpha\sigma^{-2}g\left(\left\|v_3\right\|\right)^2 + 2\alpha k_2\right)\right]^{k_2} \\
&\times \left\{1+2\alpha\sigma^{-2}g\left(\left\|v_3\right\|\right)^2\left[1+\exp\left(\alpha\sigma^{-2}g\left(\left\|v_3\right\|\right)^2 - 2\alpha k_2\right)\right]^{-1}\right\}^{-1} \\
&\times \exp\left(-(2\sigma^2)^{-1}\left\|g\left(\left\|v_3\right\|\right)v_3/\left\|v_3\right\| - (X_2'(I-P_R)X_2)^{1/2}\beta_2\right\|^2\right).
\end{aligned}$$

Since $V_1$ is independent of $V_2$, and hence of $V_3$, the joint density of $[V_1' : V_3']'$ exists and is given by

$$(2\pi\sigma^2)^{-k_1/2}\exp\{-(2\sigma^2)^{-1}\left\|v_1 - (X_1'X_1)^{-1/2}X_1'X\beta\right\|^2\}\psi(v_3).$$

Since the matrix $[C_n : D_n]$ is non-singular we obtain for the density of $\sqrt{n}(\tilde{\beta} - \beta)$

$$\begin{aligned}&(2\pi\sigma^2)^{-k_1/2} n^{-k/2} [\det(X_1'X_1)\det(X_2'(I-P_R)X_2)]^{1/2} \\ &\times \exp\left(-(2\sigma^2)^{-1}\left\| n^{-1/2}(X_1'X_1)^{1/2}(t_1 + n^{1/2}\beta_1)\right.\right. \\ &\left.\left. + n^{-1/2}(X_1'X_1)^{-1/2}X_1'X_2(t_2 + n^{1/2}\beta_2) - (X_1'X_1)^{-1/2}X_1'X\beta\right\|^2\right) \\ &\times \psi\left(n^{-1/2}(X_2'(I-P_R)X_2)^{1/2}(t_2 + n^{1/2}\beta_2)\right).\end{aligned}$$

Note that $\det(X_1'X_1)\det(X_2'(I-P_R)X_2) = \det(X'X)$. Using this, and inserting the definition of $\psi$, delivers the final result (6). □

**Remark 4.** From Proposition 1 one can immediately obtain the finite-sample distribution of $\sqrt{n}A_n(\tilde{\beta}-\beta)$ by premultiplying (3) or (4) by $A_n$. Here $A_n$ is an arbitrary (nonstochastic) $p_n \times k$ matrix. If $A_n$ has full row-rank equal to $k$ (implying $p_n = k$), this distribution has a density, which is given by $\det(A_n)^{-1} f_{n,\beta,\sigma}(A_n^{-1}s)$, $s \in \mathbb{R}^k$.

## 3. Asymptotic properties

For the asymptotic results we shall – besides the basic assumptions made in the preceding section – also assume that

$$\lim_{n\to\infty} X'X/n = Q \tag{7}$$

exists and is positive definite, i.e., $Q > 0$. We first establish "uniform $\sqrt{n}$-consistency" of the model averaging estimator, implying, in particular, uniform consistency of this estimator.

**Theorem 5.** *Suppose (7) holds.*

1. *Then $\tilde{\beta}$ is uniformly $\sqrt{n}$-consistent for $\beta$, in the sense that*

$$\lim_{M\to\infty} \sup_{n\geq k} \sup_{\beta\in\mathbb{R}^k} \mathbb{P}_{n,\beta,\sigma}\left(\sqrt{n}\left\|\tilde{\beta}-\beta\right\| \geq M\right) = 0. \tag{8}$$

   *Consequently, for every $\varepsilon > 0$*

$$\lim_{n\to\infty} \sup_{\beta\in\mathbb{R}^k} \mathbb{P}_{n,\beta,\sigma}\left(c_n\left\|\tilde{\beta}-\beta\right\| \geq \varepsilon\right) = 0 \tag{9}$$

   *holds for any sequence of real numbers $c_n \geq 0$ satisfying $c_n = o(n^{1/2})$; which reduces to uniform consistency for $c_n = 1$.*
2. *The results in Part 1 also hold for $A_n\tilde{\beta}$ as an estimator of $A_n\beta$, where $A_n$ are arbitrary (nonstochastic) matrices of dimension $p_n \times k$ such that the largest eigenvalues $\lambda_{\max}(A_n'A_n)$ are bounded.*

*Proof.* We prove (8) first. Rewrite the model averaging estimator as $\tilde{\beta} = \hat{\beta}(U) + \hat{\lambda}(\hat{\beta}(R) - \hat{\beta}(U))$. Since

$$\left\|\tilde{\beta}-\beta\right\| \leq \left\|\hat{\beta}(U)-\beta\right\| + \left|\hat{\lambda}\right| \left\|\hat{\beta}(R)-\hat{\beta}(U)\right\|,$$

since

$$\mathbb{P}_{n,\beta,\sigma}\left(\sqrt{n}\left\|\hat{\beta}(U)-\beta\right\| \geq M\right) \leq M^{-2}\sigma^2 \operatorname{trace}[(X'X/n)^{-1}],$$

and since $\operatorname{trace}[(X'X/n)^{-1}] \to \operatorname{trace}[Q^{-1}] < \infty$, it suffices to establish

$$\lim_{M\to\infty} \sup_{n\geq k} \sup_{\beta\in\mathbb{R}^k} \mathbb{P}_{n,\beta,\sigma}\left(\sqrt{n}\left|\hat{\lambda}\right| \left\|\hat{\beta}(R) - \hat{\beta}(U)\right\| \geq M\right) = 0. \tag{10}$$

Now, using (2) and the elementary inequality $z^2/[1+c\exp(z^2)]^2 \leq c^{-2}$ we have

$$\begin{aligned}
&\hat{\lambda}^2 \left\|\hat{\beta}(R) - \hat{\beta}(U)\right\|^2 \\
&\quad\leq \hat{\lambda}^2 \lambda_{\min}^{-1}(X'X) \left\|X\hat{\beta}(R) - X\hat{\beta}(U)\right\|^2 \\
&\quad= \lambda_{\min}^{-1}(X'X)\left[1+\exp(-2\alpha k_2)\exp\left(\alpha\left\|X\hat{\beta}(R) - X\hat{\beta}(U)\right\|^2/\sigma^2\right)\right]^{-2} \\
&\qquad\times \left\|X\hat{\beta}(R) - X\hat{\beta}(U)\right\|^2 \\
&\quad\leq n^{-1}\lambda_{\min}^{-1}(X'X/n)\alpha^{-1}\sigma^2\exp(4\alpha k_2) \leq Kn^{-1}\sigma^2
\end{aligned} \tag{11}$$

for a suitable finite constant $K$, since $\lambda_{\min}(X'X/n) \to \lambda_{\min}(Q) > 0$. This proves (10) and thus completes the proof of (8). The remaining claims in Part 1 follow now immediately. Part 2 is an immediate consequence of Part 1, of the inequality

$$\left\|A_n\tilde{\beta} - A_n\beta\right\|^2 \leq \lambda_{\max}(A_n'A_n)\left\|\tilde{\beta} - \beta\right\|^2,$$

and of the assumption on $\lambda_{\max}(A_n'A_n)$. □

**Remark 6.** (i) The proof has in fact shown that the difference between $\tilde{\beta}$ and $\hat{\beta}(U)$ is bounded in norm by a deterministic sequence of the form $const * \sigma n^{-1/2}$.

(ii) Although of little statistical significance since $\sigma^2$ is here assumed to be known, the proof also shows that the above proposition remains true if a supremum over $0 < \sigma^2 \leq S$, $(0 < S < \infty)$ is inserted in (8) and (9).

In the next two theorems we give the asymptotic distribution under general "moving parameter" asymptotics. Note that the case of fixed parameter asymptotics ($\beta^{(n)} \equiv \beta$) as well as the case of the usual local alternative asymptotics ($\beta^{(n)} = \beta + \delta/\sqrt{n}$) is covered by the subsequent theorems. In both these cases, Part 1 of the subsequent theorem applies if $\beta_2 \neq 0$, while Part 2 with $\gamma = 0$ and $\gamma = \delta_2$, respectively, applies if $\beta_2 = 0$.

**Theorem 7.** *Suppose (7) holds.*

1. *Let $\beta^{(n)}$ be a sequence of parameters such that $\|\sqrt{n}\beta_2^{(n)}\| \to \infty$ as $n \to \infty$. Then the distribution of $\sqrt{n}(\tilde{\beta} - \beta^{(n)})$ under $\mathbb{P}_{n,\beta^{(n)},\sigma}$ converges weakly to a $N(0, \sigma^2 Q^{-1})$-distribution.*
2. *Let $\beta^{(n)}$ be a sequence of parameters such that $\sqrt{n}\beta_2^{(n)} \to \gamma \in \mathbb{R}^{k_2}$ as $n \to \infty$. Then the distribution of $\sqrt{n}(\tilde{\beta} - \beta^{(n)})$ under $\mathbb{P}_{n,\beta^{(n)},\sigma}$ converges weakly to the distribution of*

$$\begin{aligned}
&B_\infty\gamma + C_\infty Z_1 \\
&+ \left[1+\exp(2\alpha k_2)\exp\left(-\alpha\left\|Z_2 + (Q_{22} - Q_{21}Q_{11}^{-1}Q_{12})^{1/2}\gamma\right\|^2/\sigma^2\right)\right]^{-1} \\
&\times\{D_\infty Z_2 - B_\infty\gamma\}
\end{aligned} \tag{12}$$

*where*

$$B_\infty = \begin{bmatrix} Q_{11}^{-1}Q_{12} \\ -I_{k_2} \end{bmatrix}, \qquad C_\infty = \begin{bmatrix} Q_{11}^{-1/2} \\ 0_{k_2\times k_1} \end{bmatrix},$$
$$D_\infty = \begin{bmatrix} -Q_{11}^{-1}Q_{12}(Q_{22}-Q_{21}Q_{11}^{-1}Q_{12})^{-1/2} \\ (Q_{22}-Q_{21}Q_{11}^{-1}Q_{12})^{-1/2} \end{bmatrix},$$

*and where* $Z_1 \sim N(0,\sigma^2 I_{k_1})$ *is independent of* $Z_2 \sim N(0,\sigma^2 I_{k_2})$. *The density of the distribution of (12) is given by*

$$\begin{aligned} f_{\infty,\gamma}(t) = {} & (2\pi\sigma^2)^{-k/2}[\det(Q)]^{1/2} \\ & \times \exp\left(-(2\sigma^2)^{-1}\left\|Q_{11}^{1/2}t_1 + Q_{11}^{-1/2}Q_{12}t_2\right\|^2\right) \\ & \times \left[1+\exp\left(-\alpha\sigma^{-2}g\left(\left\|D_{\infty 2}^{-1}(t_2+\gamma)\right\|\right)^2 + 2\alpha k_2\right)\right]^{k_2} \\ & \times \left\{1+2\alpha\sigma^{-2}g\left(\left\|D_{\infty 2}^{-1}(t_2+\gamma)\right\|\right)^2\right. \\ & \times \left.\left[1+\exp\left(\alpha\sigma^{-2}g\left(\left\|D_{\infty 2}^{-1}(t_2+\gamma)\right\|\right)^2 - 2\alpha k_2\right)\right]^{-1}\right\}^{-1} \\ & \times \exp\left\{-(2\sigma^2)^{-1}\left\|g\left(\left\|D_{\infty 2}^{-1}(t_2+\gamma)\right\|\right)\left\|D_{\infty 2}^{-1}(t_2+\gamma)\right\|^{-1}\right.\right. \\ & \times \left.\left. D_{\infty 2}^{-1}(t_2+\gamma) - D_{\infty 2}^{-1}\gamma\right\|^2\right\}, \end{aligned} \tag{13}$$

*where $t$ is partitioned as $(t_1', t_2')'$ with $t_1$ being a $k_1 \times 1$ vector. Furthermore,* $D_{\infty 2} = (Q_{22} - Q_{21}Q_{11}^{-1}Q_{12})^{-1/2}$, *and $g$ is as defined in the Appendix (with* $a = \exp(2\alpha k_2)$ *and* $b = \alpha^{-1}\sigma^2$*).*

*Proof.* To prove Part 1 represent $\sqrt{n}(\tilde\beta - \beta^{(n)})$ as $\sqrt{n}(\hat\beta(U) - \beta^{(n)}) + \hat\lambda\sqrt{n}(\hat\beta(R) - \hat\beta(U))$. The first term is $N(0,\sigma^2(X'X/n)^{-1})$-distributed under $\mathbb{P}_{n,\beta^{(n)},\sigma}$, which obviously converges to a $N(0,\sigma^2 Q^{-1})$-distribution. It hence suffices to show that $\hat\lambda\sqrt{n}(\hat\beta(R) - \hat\beta(U))$ converges to zero in $\mathbb{P}_{n,\beta^{(n)},\sigma}$-probability. Since $\lambda_{\min}^{-1}(X'X/n)$ is bounded by assumption (7) and since

$$\begin{aligned} \hat\lambda^2\left\|\sqrt{n}(\hat\beta(R) - \hat\beta(U))\right\|^2 \le {} & n\lambda_{\min}^{-1}(X'X)\left\|X\hat\beta(R) - X\hat\beta(U)\right\|^2 \\ & \times\left[1+\exp\left(\alpha\sigma^{-2}\left\|X\hat\beta(R) - X\hat\beta(U)\right\|^2 - 2\alpha k_2\right)\right]^{-2} \end{aligned}$$

as shown in (11), it furthermore suffices to show that

$$\left\|X\hat\beta(R) - X\hat\beta(U)\right\|^2 \to \infty \text{ in } \mathbb{P}_{n,\beta^{(n)},\sigma}\text{-probability.} \tag{14}$$

Note that

$$\begin{aligned} \left\|X\hat\beta(R) - X\hat\beta(U)\right\|^2 &= \|(P_U - P_R)Y\|^2 \\ &= \left\|(P_U - P_R)u + (P_U - P_R)X_2\beta_2^{(n)}\right\|^2 \\ &\ge \left|\left\|(P_U - P_R)X_2\beta_2^{(n)}\right\| - \|(P_U - P_R)u\|\right|^2. \end{aligned}$$

The second term satisfies $\mathbb{E}_{n,\beta^{(n)},\sigma}\|(P_U - P_R)u\|^2 = \sigma^2 k_2$ and hence is stochastically bounded in $\mathbb{P}_{n,\beta^{(n)},\sigma}$-probability. The square of the first term, i.e.,

$$\left\|(P_U - P_R)X_2\beta_2^{(n)}\right\|^2$$

equals

$$\sqrt{n}\beta_2^{(n)\prime}[(X_2'X_2/n) - (X_2'X_1/n)(X_1'X_1/n)^{-1}(X_1'X_2/n)]\sqrt{n}\beta_2^{(n)}.$$

Since the matrix in brackets converges to $Q_{22} - Q_{21}Q_{11}Q_{12}$, which is positive definite, the above display diverges to infinity, establishing (14). This completes the proof of Part 1.

We next turn to the proof of Part 2. The proof of (12) is immediate from (3) upon observing that $B_n \to B_\infty$, $\sqrt{n}C_n \to C_\infty$, and $\sqrt{n}D_n \to D_\infty$. To prove (13) observe that (12) can be written as

$$\begin{aligned}
&B_\infty\gamma + C_\infty Z_1 \\
&+ \left[1 + \exp(2\alpha k_2)\exp\left(-\alpha\left\|Z_2 + (Q_{22} - Q_{21}Q_{11}^{-1}Q_{12})^{1/2}\gamma\right\|^2/\sigma^2\right)\right]^{-1} \\
&\times \{D_\infty(Z_2 + (Q_{22} - Q_{21}Q_{11}^{-1}Q_{12})^{1/2}\gamma)\} \\
&\qquad = B_\infty\gamma + C_\infty Z_1 + D_\infty\left[1 + \exp(2\alpha k_2)\exp\left(-\alpha\|W_2\|^2/\sigma^2\right)\right]^{-1} W_2
\end{aligned}$$

where $W_2 \sim N((Q_{22} - Q_{21}Q_{11}^{-1}Q_{12})^{1/2}\gamma, \sigma^2 I_{k_2})$ is independent of $Z_1$. Again using Lemmata 15 and 16 in the Appendix gives the density of

$$W_3 = \left[1 + \exp(2\alpha k_2)\exp\left(-\alpha\|W_2\|^2/\sigma^2\right)\right]^{-1} W_2$$

as

$$\begin{aligned}
\chi(w_3) = {} & (2\pi\sigma^2)^{-k_2/2}\left[1 + \exp\left(-\alpha\sigma^{-2}g(\|w_3\|)^2 + 2\alpha k_2\right)\right]^{k_2} \\
& \times \left\{1 + 2\alpha\sigma^{-2}g(\|w_3\|)^2\left[1 + \exp\left(\alpha\sigma^{-2}g(\|w_3\|)^2 - 2\alpha k_2\right)\right]^{-1}\right\}^{-1} \\
& \times \exp\left(-(2\sigma^2)^{-1}\left\|g(\|w_3\|)w_3/\|w_3\| - (Q_{22} - Q_{21}Q_{11}^{-1}Q_{12})^{1/2}\gamma\right\|^2\right).
\end{aligned}$$

Since $Z_1$ is independent of $Z_2$, and hence of $W_3$, the joint density of $[Z_1' : W_3']'$ exists and is given by

$$(2\pi\sigma^2)^{-k_1/2}\exp\left(-(2\sigma^2)^{-1}\|z_1\|^2\right)\chi(w_3).$$

Since the matrix $[C_\infty : D_\infty]$ is non-singular we obtain finally

$$\begin{aligned}
& (2\pi\sigma^2)^{-k_1/2}\left[\det(Q_{11})\det(Q_{22} - Q_{21}Q_{11}^{-1}Q_{12})\right]^{1/2} \\
& \times \exp\left(-(2\sigma^2)^{-1}\left\|Q_{11}^{1/2}(t_1 - Q_{11}^{-1}Q_{12}\gamma) + Q_{11}^{-1/2}Q_{12}(t_2 + \gamma)\right\|^2\right) \\
& \times \chi\left((Q_{22} - Q_{21}Q_{11}^{-1}Q_{12})^{1/2}(t_2 + \gamma)\right).
\end{aligned}$$

Inserting the expression for $\chi$ derived above gives (13). □

Since in both cases considered in the above theorem the limiting distribution is continuous, the finite-sample cumulative distribution function (cdf)

$$F_{n,\beta^{(n)},\sigma}(t) = \mathbb{P}_{n,\beta^{(n)},\sigma}\left(\sqrt{n}(\tilde{\beta} - \beta^{(n)}) \leq t\right)$$

converges to the cdf of the corresponding limiting distribution even in the sup-norm as a consequence of the multivariate version of Polya's Theorem (cf. [1], Ex.6, [2]). We next show that the convergence occurs in an even stronger sense. Let $f_\infty$ denote the density of the asymptotic distribution of $\sqrt{n}(\tilde{\beta} - \beta^{(n)})$ given in the previous theorem. That is, $f_\infty$ is equal to $f_{\infty,\gamma}$ given in (13) if $\sqrt{n}\beta_2^{(n)} \to \gamma \in \mathbb{R}^{k_2}$, and is equal to the density of an $N(0, \sigma^2 Q^{-1})$-distribution if $\|\sqrt{n}\beta_2^{(n)}\| \to \infty$. For obvious reasons and for convenience we shall denote the $N(0, \sigma^2 Q^{-1})$-density by $f_{\infty,\infty}$.

**Theorem 8.** *Suppose the assumptions of Theorem 7 hold. Then the finite-sample density $f_{n,\beta^{(n)},\sigma}$ of $\sqrt{n}(\tilde{\beta} - \beta^{(n)})$ converges to $f_\infty$, the density of the corresponding asymptotic distribution, in the $L^1$-sense. Consequently, the finite-sample cdf $F_{n,\beta^{(n)},\sigma}$ converges to the corresponding asymptotic cdf in total variation distance.*

*Proof.* In the case where $\sqrt{n}\beta_2^{(n)} \to \gamma \in \mathbb{R}^{k_2}$, inspection of (6), and noting that $g$ as well as $T^{-1}$ given in Lemma 15 are continuous, shows that (6) converges to (13) pointwise. In the case where $\|\sqrt{n}\beta_2^{(n)}\| \to \infty$, Lemma 17 in the Appendix and inspection of (6) show that (6) converges pointwise to the density of a $N(0, \sigma^2 Q^{-1})$-distribution. Observing that $f_{n,\beta^{(n)},\sigma}$ as well as $f_\infty$ are probability densities, the proof is then completed by an application of Scheffé's lemma. □

**Remark 9.** We note for later use that inspection of (13) combined with Lemma 17 in the Appendix shows that for $\|\gamma\| \to \infty$ we have $f_{\infty,\gamma} \to f_{\infty,\infty}$ (the $N(0, \sigma^2 Q^{-1})$-density) pointwise on $\mathbb{R}^k$, and hence also in the $L^1$-sense. As a consequence, the corresponding cdfs converge in the total variation sense to the cdf of a $N(0, \sigma^2 Q^{-1})$-distribution.

**Remark 10.** The results in this section imply that the convergence of the finite-sample cdf to the asymptotic cdf does not occur uniformly w.r.t. the parameter $\beta$. [Cf. also the first step in the proof of Theorem 13 below.]

**Remark 11.** Theorems 7 and 8 in fact provide a characterization of all accumulation points of the finite sample distribution $F_{n,\beta^{(n)},\sigma}$ (w.r.t. the total variation topology) for arbitrary sequences $\beta^{(n)}$. This follows from a simple subsequence argument applied to $\sqrt{n}\beta_2^{(n)}$ and observing that $(\mathbb{R} \cup \{-\infty, \infty\})^{k_2}$ is compact; cf. also Remark 4.4 in [7].

**Remark 12.** Part 1 of Theorem 7 as well as the representation (12) immediately generalize to $\sqrt{n}A(\tilde{\beta} - \beta)$ with $A$ a non-stochastic $p \times k$ matrix. If $A$ has full row-rank equal to $k$, the resulting asymptotic distribution has a density, which is given by $\det(A)^{-1} f_\infty(A^{-1}s)$, $s \in \mathbb{R}^k$.

## 4. Estimation of the finite-sample distribution: an impossibility result

As can be seen from Theorem 3, the finite-sample distribution depends on the unknown parameter $\beta$, even after centering at $\beta$. Hence, it is obviously of interest to estimate this distribution, e.g., for purposes of conducting inference. It is easy to construct a consistent estimator of the cumulative distribution function $F_{n,\beta,\sigma}$ of the scaled and centered model averaging estimator $\tilde{\beta}$, i.e., of

$$F_{n,\beta,\sigma}(t) = \mathbb{P}_{n,\beta,\sigma}\left(\sqrt{n}(\tilde{\beta}-\beta) \leq t\right).$$

To this end, let $\hat{M}$ be an estimator that consistently decides between the restricted model $M_R$ and the unrestricted model $M_U$, i.e., $\lim_{n\to\infty} \mathbb{P}_{n,\beta,\sigma}(\hat{M} = M_R) = 1$ if $\beta_2 = 0$ and $\lim_{n\to\infty} \mathbb{P}_{n,\beta,\sigma}(\hat{M} = M_U) = 1$ if $\beta_2 \neq 0$. [Such a procedure is easily constructed, e.g., from BIC or from a $t$-test for the hypothesis $\beta_2 = 0$ with critical value that diverges to infinity at a rate slower than $n^{1/2}$.] Define $\check{f}_n$ equal to $f^{\dagger}_{\infty,\infty}$, the density of the $N(0, \sigma^2(X'X/n)^{-1})$-distribution, on the event $\hat{M} = M_U$, and define $\check{f}_n$ equal to $f^{\dagger}_{\infty,0}$ otherwise, where $f^{\dagger}_{\infty,0}$ follows the same formula as $f_{\infty,0}$, with the only exception that $Q$ is replaced by $X'X/n$. Then – as is proved in the Appendix –

$$\int_{\mathbb{R}^k} \left|\check{f}_n(z) - f_{n,\beta,\sigma}(z)\right| dz \to 0 \tag{15}$$

in $\mathbb{P}_{n,\beta,\sigma}$-probability as $n \to \infty$ for every $\beta \in \mathbb{R}^k$. Define $\check{F}_n$ as the cdf corresponding to $\check{f}_n$. Then for every $\delta > 0$

$$\mathbb{P}_{n,\beta,\sigma}\left(\left\|\check{F}_n - F_{n,\beta,\sigma}\right\|_{TV} > \delta\right) \to 0$$

as $n \to \infty$, where $\|\cdot\|_{TV}$ denotes the total variation norm. This shows that $\check{F}_n$ is a consistent estimator of $F_{n,\beta,\sigma}$ in the total variation distance. A fortiori then also

$$\mathbb{P}_{n,\beta,\sigma}\left(\sup_t \left|\check{F}_n(t) - F_{n,\beta,\sigma}(t)\right| > \delta\right) \to 0$$

holds.

The estimator $\check{F}_n$ just constructed has been obtained from the asymptotic cdf by replacing unknown quantities with suitable estimators. As noted in Remark 10, the convergence of the finite-sample cdf to their asymptotic counterpart does not occur uniformly w.r.t. the parameter $\beta$. Hence, it is to be expected that $\check{F}_n$ will inherit this deficiency, i.e., $\check{F}_n$ will not be uniformly consistent. Of course, this makes it problematic to base inference on $\check{F}_n$, as then there is no guarantee – at *any* sample size – that $\check{F}_n$ will be close to the true cdf. This naturally raises the question if estimators other than $\check{F}_n$ exist that are uniformly consistent. The answer turns out to be negative as we show in the next theorem. In fact, uniform consistency fails dramatically, cf. (17) below. This result further shows that uniform consistency already fails over certain shrinking balls in the parameter space (and thus a fortiori fails in general over compact subsets of the parameter space), and fails even if one considers the easier estimation problem of estimating $F_{n,\beta,\sigma}$ only at a given value of the argument $t$ rather than estimating the entire function $F_{n,\beta,\sigma}$ (and measuring loss in a norm like the total variation norm or the sup-norm). Although of little statistical significance, we note that a similar result can be obtained for the problem of estimating the asymptotic cdf. Related impossibility results for post-model-selection estimators as well as for certain shrinkage-type estimators are given in [9, 10, 11].

In the result to follow we shall consider estimators of $F_{n,\beta,\sigma}(t)$ at a *fixed* value of the argument $t$. An estimator of $F_{n,\beta,\sigma}(t)$ is now nothing else than a real-valued random variable $\Gamma_n = \Gamma_n(Y, X)$. For mnemonic reasons we shall, however, use the symbol $\hat{F}_n(t)$ instead of $\Gamma_n$ to denote an arbitrary estimator of $F_{n,\beta,\sigma}(t)$. This notation should not be taken as implying that the estimator is obtained by evaluating an estimated cdf at the argument $t$, or that it is constrained to lie between zero and one. For simplicity, we give the impossibility result only in the simple situation where $k_2 = 1$ and $Q$ is block-diagonal, i.e., $X_1$ and $X_2$ are asymptotically orthogonal. There is no reason to believe that the non-uniformity problem will disappear in more complicated situations.

**Theorem 13.** *Suppose (7) holds. Suppose further that $k_2 = 1$ and that $Q$ is block-diagonal, i.e., the $k_1 \times k_2$ matrix $Q_{12}$ is equal to zero. Then the following holds for every $\beta \in M_R$ and every $t \in \mathbb{R}^k$: There exist $\delta_0 > 0$ and $\rho_0$, $0 < \rho_0 < \infty$, such that any estimator $\hat{F}_n(t)$ of $F_{n,\beta,\sigma}(t)$ satisfying*

$$\mathbb{P}_{n,\beta,\sigma}\left(\left|\hat{F}_n(t) - F_{n,\beta,\sigma}(t)\right| > \delta\right) \overset{n\to\infty}{\longrightarrow} 0 \tag{16}$$

*for every $\delta > 0$ (in particular, every estimator that is consistent) also satisfies*

$$\sup_{\substack{\vartheta\in\mathbb{R}^k \\ ||\vartheta-\beta||<\rho_0/\sqrt{n}}} \mathbb{P}_{n,\vartheta,\sigma}\left(\left|\hat{F}_n(t) - F_{n,\vartheta,\sigma}(t)\right| > \delta_0\right) \overset{n\to\infty}{\longrightarrow} 1. \tag{17}$$

*The constants $\delta_0$ and $\rho_0$ may be chosen in such a way that they depend only on $t$, $Q$, $\sigma$, and the tuning parameter $\alpha$. Moreover,*

$$\liminf_{n\to\infty} \inf_{\hat{F}_n(t)} \sup_{\substack{\vartheta\in\mathbb{R}^k \\ ||\vartheta-\beta||<\rho_0/\sqrt{n}}} \mathbb{P}_{n,\vartheta,\sigma}\left(\left|\hat{F}_n(t) - F_{n,\vartheta,\sigma}(t)\right| > \delta_0\right) > 0 \tag{18}$$

*and*

$$\sup_{\delta>0} \liminf_{n\to\infty} \inf_{\hat{F}_n(t)} \sup_{\substack{\vartheta\in\mathbb{R}^k \\ ||\vartheta-\beta||<\rho_0/\sqrt{n}}} \mathbb{P}_{n,\vartheta,\sigma}\left(\left|\hat{F}_n(t) - F_{n,\vartheta,\sigma}(t)\right| > \delta\right) \geq \frac{1}{2}, \tag{19}$$

*where the infima in (18) and (19) extend over* all *estimators $\hat{F}_n(t)$ of $F_{n,\beta,\sigma}(t)$.*

*Proof. Step 1*: Let $\beta \in M_R$ and $t \in \mathbb{R}^k$ be given. Observe that by Theorems 7 and 8 the limit

$$F_{\infty,\gamma}(t) := \lim F_{n,\beta+(\eta,\gamma)'/\sqrt{n},\sigma}(t)$$

exists for every $\eta \in \mathbb{R}^{k_1}$, $\gamma \in \mathbb{R}^{k_2} = \mathbb{R}$, and does not depend on $\eta$. We now show that $F_{\infty,\gamma}(t)$ is non-constant in $\gamma \in \mathbb{R}$. First, observe that by Remark 9 and the block-diagonality assumption on $Q$

$$\lim_{\|\gamma\|\to\infty} F_{\infty,\gamma}(t) = \mathbb{P}\left(Q_{11}^{-1/2} Z_1 \leq t_1\right) \mathbb{P}\left(Q_{22}^{-1/2} Z_2 \leq t_2\right)$$

where $Z_1$ and $Z_2$ are as in Theorem 7, $t$ is partitioned as $(t_1', t_2)'$ with $t_2$ a scalar, and $\mathbb{P}$ is the probability measure governing $(Z_1', Z_2)'$. Second, we have from (12) and the block-diagonality assumption on $Q$ that $F_{\infty,\gamma}(t)$ is the product of

$$\mathbb{P}\left(Q_{11}^{-1/2} Z_1 \leq t_1\right)$$

with

$$\begin{aligned} \mathbb{P}\Bigg(&\left[1 + \exp(2\alpha)\exp\left(-\alpha\left(Z_2 + Q_{22}^{1/2}\gamma\right)^2/\sigma^2\right)\right]^{-1} \\ &\times \left(Q_{22}^{-1/2} Z_2 + \gamma\right) - \gamma \leq t_2\Bigg). \end{aligned} \tag{20}$$

Since $\mathbb{P}(Q_{11}^{-1/2}Z_1 \leq t_1)$ is positive and independent of $\gamma$, it suffices to show that (20) differs from $\mathbb{P}(Q_{22}^{-1/2}Z_2 \leq t_2)$ for at least one $\gamma \in \mathbf{R}$. Suppose first that $t_2 > 0$. Then specializing to the case $\gamma = 0$ in (20) it suffices to show that

$$\mathbb{P}\left(\left[1+\exp(2\alpha)\exp\left(-\alpha Z_2^2/\sigma^2\right)\right]^{-1} Q_{22}^{-1/2} Z_2 \leq t_2\right). \tag{21}$$

differs from $\mathbb{P}(Q_{22}^{-1/2}Z_2 \leq t_2)$. But this follows from

$$\begin{aligned}
&\mathbb{P}\left(\left[1+\exp(2\alpha)\exp\left(-\alpha Z_2^2/\sigma^2\right)\right]^{-1} Q_{22}^{-1/2} Z_2 \leq t_2\right)\\
&\quad = 1/2 + \mathbb{P}\left(Z_2 \geq 0, h(Z_2) \leq Q_{22}^{1/2} t_2\right)\\
&\quad = 1/2 + \mathbb{P}\left(0 \leq Z_2 \leq g\left(Q_{22}^{1/2} t_2\right)\right)\\
&\quad > 1/2 + \mathbb{P}\left(0 \leq Z_2 \leq Q_{22}^{1/2} t_2\right)\\
&\quad = \mathbb{P}\left(Q_{22}^{-1/2} Z_2 \leq t_2\right)
\end{aligned}$$

since $h$ as defined in the Appendix (with $a = \exp(2\alpha)$ and $b = \sigma^2/\alpha$) is strictly monotonically increasing and satisfies $h(x) < x$ for every $x > 0$, which entails $g(y) > y$ for every $y > 0$. For symmetry reasons a dual statement holds for $t_2 < 0$. It remains to consider the case $t_2 = 0$. In this case (20) equals

$$\begin{aligned}
&\mathbb{P}\left(\left[1+\exp(2\alpha)\exp\left(-\alpha\left(Z_2 + Q_{22}^{1/2}\gamma\right)^2/\sigma^2\right)\right]^{-1}\right.\\
&\qquad \left. \times \left(Z_2 + Q_{22}^{1/2}\gamma\right) \leq Q_{22}^{1/2}\gamma\right).
\end{aligned} \tag{22}$$

Let $\gamma > 0$ be arbitrary. Then (22) equals

$$\mathbb{P}\left(Z_2 + Q_{22}^{1/2}\gamma < 0\right) + \mathbb{P}\left(Z_2 + Q_{22}^{1/2}\gamma \geq 0, h\left(Z_2 + Q_{22}^{1/2}\gamma\right) \leq Q_{22}^{1/2}\gamma\right).$$

Arguing as before, this can be written as

$$\begin{aligned}
&\mathbb{P}\left(Z_2 + Q_{22}^{1/2}\gamma < 0\right) + \mathbb{P}\left(0 \leq Z_2 + Q_{22}^{1/2}\gamma \leq g\left(Q_{22}^{1/2}\gamma\right)\right)\\
&\quad > \mathbb{P}\left(Z_2 + Q_{22}^{1/2}\gamma < 0\right) + \mathbb{P}\left(0 \leq Z_2 + Q_{22}^{1/2}\gamma \leq Q_{22}^{1/2}\gamma\right)\\
&\quad = \mathbb{P}\left(Z_2 \leq 0\right) = \mathbb{P}\left(Q_{22}^{-1/2} Z_2 \leq 0\right)
\end{aligned}$$

which completes the proof of Step 1.

*Step 2*: We prove (17) and (18) first. For this purpose we make use of Lemma 3.1 in Leeb and Pötscher [11] with the notational identification $\alpha = \beta \in M_R$, $B = \mathbb{R}^k$, $B_n = \{\vartheta \in \mathbb{R}^k : \|\vartheta - \beta\| < \rho_0 n^{-1/2}\}$, $\varphi_n(\cdot) = F_{n,\cdot,\sigma}(t)$, and $\hat{\varphi}_n = \hat{F}_n(t)$, where $\rho_0$ will be chosen shortly. The contiguity assumption of this lemma is obviously satisfied; cf. also Lemma A.1 in [11]. It hence remains to show that there exists a value of $\rho_0$, $0 < \rho_0 < \infty$, such that $\delta^*$ defined in Lemma 3.1 of Leeb and Pötscher [11], which represents the limit inferior of the oscillation of $F_{n,\cdot,\sigma}(t)$ over $B_n$, is positive. Applying Lemma 3.5(a) of Leeb and Pötscher [11] with $\zeta_n = \rho_0 n^{-1/2}$ and the set $G_0$ equal to $G = \{(\eta', \gamma)' \in \mathbb{R}^k : \|(\eta', \gamma)'\| < 1\}$, it suffices to show that $F_{\infty,\gamma}(t)$ viewed as a function of $(\eta', \gamma)'$ is non-constant on the set $\{(\eta', \gamma)' \in \mathbb{R}^k :$

$\|(\eta', \gamma)'\| < \rho_0\}$; in view of Lemma 3.1 of Leeb and Pötscher [11], the corresponding $\delta_0$ can then be chosen as any positive number less than one-half of the oscillation of $F_{\infty,\gamma}(t)$ over this set. That such a $\rho_0$ indeed exists now follows from Step 1. Furthermore, observe that $F_{\infty,\cdot}(t)$ depends only on $\alpha$, $Q$, $\sigma$, and $t$. Hence, $\delta_0$ and $\rho_0$ may be chosen such that they also only depend on these quantities. This completes the proof of (17) and (18).

To prove (19) we use Corollary 3.4 in [11] with the same identification of notation as above, with $\zeta_n = \rho_0 n^{-1/2}$, and with $V = \mathbb{R}^k$. The asymptotic uniform equicontinuity condition in that corollary is then satisfied in view of

$$\|\mathbb{P}_{n,\theta,\sigma} - \mathbb{P}_{n,\vartheta,\sigma}\|_{TV} \leq 2\Phi\left(\|\theta - \vartheta\| \lambda_{\max}^{1/2}(X'X)/(2\sigma)\right) - 1,$$

cf. Lemma A.1 in [11]. Given that the positivity of $\delta^*$ has already been established in the previous paragraph, applying Corollary 3.4 in [11] then establishes (19). □

**Remark 14.** The impossibility result given in the above theorem also holds for the class of randomized estimators (with $\mathbb{P}_{n,\cdot,\sigma}$ replaced by $\mathbb{P}^*_{n,\cdot,\sigma}$, the distribution of the randomized sample). This follows immediately from Lemma 3.6 in [11] and the attending discussion.

## Appendix A: Some technical results

Let the function $h : [0, \infty) \to [0, \infty)$ be given by $h(\xi) = [1 + a\exp(-\xi^2/b)]^{-1}\xi$ where $a$ and $b$ are positive real numbers. It is easy to see that $h$ is strictly monotonically increasing on $[0, \infty)$, is continuous, satisfies $h(0) = 0$ and $\lim_{\xi\to\infty} h(\xi) = \infty$. The inverse $g : [0, \infty) \to [0, \infty)$ of $h$ clearly exists, is strictly monotonically increasing on $[0, \infty)$, is continuous, satisfies $g(0) = 0$ and $\lim_{\zeta\to\infty} g(\zeta) = \infty$. In the following lemma we shall use the natural convention that $g(\|y\|)y/\|y\| = 0$ for $y = 0$, which makes $y \to g(\|y\|)y/\|y\|$ a continuous function on all of $\mathbb{R}^m$.

**Lemma 15.** *Let $T : \mathbb{R}^m \to \mathbb{R}^m$ be given by*

$$T(x) = \left[1 + a\exp(-\|x\|^2/b)\right]^{-1} x$$

*where $a$ and $b$ are positive real numbers. Then $T$ is a bijection. Its inverse is given by*

$$T^{-1}(y) = g(\|y\|)y/\|y\|$$

*where $g$ has been defined above. Moreover, $T^{-1}$ is continuously partially differentiable and $\|T^{-1}(y)\| = g(\|y\|)$ holds for all $y$.*

*Proof.* If $y = 0$ it is obvious that $T(T^{-1}(y)) = 0 = y$ in view of the convention made above. Now suppose that $y \neq 0$. Then

$$\begin{aligned} T(T^{-1}(y)) &= [1 + a\exp\left(-g(\|y\|)^2/b\right)]^{-1} g(\|y\|)y/\|y\| \\ &= h\left(g(\|y\|)\right) y/\|y\| = y. \end{aligned}$$

Similarly, if $x = 0$ then $T^{-1}(T(x)) = 0$. Now suppose $x \neq 0$. Then $T(x) \neq 0$ and, observing that $\|T(x)\| = [1 + a\exp(-\|x\|^2/b)]^{-1}\|x\|$, we have

$$\begin{aligned} T^{-1}(T(x)) &= g\left(\|T(x)\|\right) T(x)/\|T(x)\| \\ &= g\left(\left[1 + a\exp\left(-\|x\|^2/b\right)\right]^{-1}\|x\|\right) x/\|x\| \\ &= g\left(h(\|x\|)\right) x/\|x\| = x. \end{aligned}$$

That $T^{-1}$ is continuously partially differentiable follows from the corresponding property of $T$ and the fact that the determinant of the derivative of $T$ does never vanish as shown in the next lemma. The final claim is obvious in case $y \neq 0$, and follows from the convention made above and the fact that $g(0) = 0$ in case $y = 0$. □

**Lemma 16.** *Let $T$ be as in the preceding lemma. Then the determinant of the derivative $D_x T$ is given by*

$$\left[1 + a\exp\left(-\|x\|^2/b\right)\right]^{-m}\left\{1 + 2b^{-1}\left[1 + a^{-1}\exp\left(\|x\|^2/b\right)\right]^{-1}\|x\|^2\right\}$$

*which is always positive.*

*Proof.* Elementary calculations show that

$$\begin{aligned} D_x T = {}& \left[1 + a\exp\left(-\|x\|^2/b\right)\right]^{-1} \\ & \times\left\{I_m + 2ab^{-1}\exp\left(-\|x\|^2/b\right)\left[1 + a\exp\left(-\|x\|^2/b\right)\right]^{-1} xx'\right\}. \end{aligned}$$

Since the determinate of $I_m + cxx'$ equals $1 + cx'x$, the result follows. □

**Lemma 17.** *For $g$ defined above we have*

$$\lim_{\zeta\to\infty} g(\zeta)/\zeta = 1$$

*and*

$$\lim_{\zeta\to\infty} \left((g(\zeta)/\zeta) - 1\right)\zeta = 0.$$

*Proof.* It suffices to prove the second claim:

$$\begin{aligned} \lim_{\zeta\to\infty} \left((g(\zeta)/\zeta) - 1\right)\zeta &= \lim_{\zeta\to\infty} (g(\zeta) - \zeta) = \lim_{\xi\to\infty} (g(h(\xi)) - h(\xi)) \\ &= \lim_{\xi\to\infty} \left(\xi - \left[1 + a\exp\left(-\xi^2/b\right)\right]^{-1}\xi\right) \\ &= \lim_{\xi\to\infty} \xi\left[1 + a^{-1}\exp\left(\xi^2/b\right)\right]^{-1} = 0. \end{aligned}$$

□

*Proof* (Verification of (15) in Section 5). In view of Theorem 8 it suffices to show that

$$\int_{\mathbb{R}^k} \left|\check{f}_n(z) - f_\infty(z)\right| dz \to 0$$

in $\mathbb{P}_{n,\beta,\sigma}$-probability as $n \to \infty$ for every $\beta \in \mathbb{R}^k$ where we recall that $f_\infty$ is equal to $f_{\infty,\infty}$, the density of an $N(0, \sigma^2 Q^{-1})$-distribution, if $\beta_2 \neq 0$, and is equal to $f_{\infty,0}$ given in (13) if $\beta_2 = 0$. Now,

$$\begin{aligned}
&\mathbb{P}_{n,\beta,\sigma}\left(\int_{\mathbb{R}^k}\left|\check{f}_n(z)-f_\infty(z)\right|dz>\varepsilon\right)\\
&\quad=\mathbb{P}_{n,\beta,\sigma}\left(\int_{\mathbb{R}^k}\left|\check{f}_n(z)-f_\infty(z)\right|dz>\varepsilon,\hat{M}=M_R\right)\\
&\qquad+\mathbb{P}_{n,\beta,\sigma}\left(\int_{\mathbb{R}^k}\left|\check{f}_n(z)-f_\infty(z)\right|dz>\varepsilon,\hat{M}=M_U\right)\\
&\quad=\mathbb{P}_{n,\beta,\sigma}\left(\int_{\mathbb{R}^k}\left|f^{\dagger}_{\infty,0}(z)-f_\infty(z)\right|dz>\varepsilon,\hat{M}=M_R\right)\\
&\qquad+\mathbb{P}_{n,\beta,\sigma}\left(\int_{\mathbb{R}^k}\left|f^{\dagger}_{\infty,\infty}(z)-f_\infty(z)\right|dz>\varepsilon,\hat{M}=M_U\right)
\end{aligned}$$

where we have made use of the definition of $\check{f}_n$. If $\beta \in M_R$, then clearly the event $\hat{M} = M_U$ has probability approaching zero and hence the last probability in the above display converges to zero. Furthermore, if $\beta \in M_R$, the last but one probability reduces to

$$\mathbb{P}_{n,\beta,\sigma}\left(\int_{\mathbb{R}^k}\left|f^{\dagger}_{\infty,0}(z)-f_{\infty,0}(z)\right|dz>\varepsilon,\hat{M}=M_R\right)$$

which converges to zero since

$$\int_{\mathbb{R}^k}\left|f^{\dagger}_{\infty,0}(z)-f_{\infty,0}(z)\right|dz\to 0$$

in view of pointwise convergence of $f^{\dagger}_{\infty,0}$ to $f_{\infty,0}$ and Scheffé's lemma. [To be able to apply Scheffé's lemma we need to know that not only $f_{\infty,0}$ but also $f^{\dagger}_{\infty,0}(z)$ is a probability density. But this is obvious, as (13) defines a probability density for *any* symmetric and positive definite matrix $Q$.] The proof for the case where $\beta \in M_U$ is completely analogous noting that then $f_\infty = f_{\infty,\infty}$ holds. □

## Acknowledgments

I would like to thank Hannes Leeb, Richard Nickl, and two anonymous referees for helpful comments on the paper.